\documentclass {article}

\usepackage {amsmath}
\usepackage {amsfonts}
\usepackage {amssymb}
\usepackage {amsthm}
\usepackage [all,cmtip] {xy}

\newtheorem {thm}   		{Theorem}
\newtheorem {prop} 		{Proposition}
\newtheorem {lem} 		{Lemma}
\newtheorem {cor} 		{Corollary}

\newtheorem {conj}              {Conjecture}

\theoremstyle {definition}
\newtheorem {define}    	{Definition}

\theoremstyle {remark}

\newcommand {\bb} 		{\mathbb}

\newcommand {\fr} 		{\mathfrak}

\newcommand {\al} 		{\alpha}

\newcommand {\dde} 		{\Delta}

\newcommand {\si} 		{\sigma}

\newcommand {\p} 			{\prime}
\newcommand {\vp} 		{\varphi}

\newcommand {\eps} 		{\epsilon}

\newcommand {\inject} {\hookrightarrow}

\newcommand {\wt}			{\widetilde}

\newcommand {\surject}{\twoheadrightarrow}

\newcommand {\Hom} 		{\text{Hom}}

\newcommand {\spec} 	{\text{Spec }}

\newcommand {\card} 	{\text{Card}}

\newcommand {\ext}    {\text{Ext}}
\newcommand {\en}		  {\text{End}}

\usepackage {amsmath}
\usepackage {amsfonts}
\usepackage {amssymb}
\usepackage {amsthm}
\usepackage [all,cmtip] {xy}

\newcommand {\add} {\mathbb{G}_a}
\newcommand {\padd} {\mathbb{G}_{a,\text{p}}}

\newcommand {\adds} {\mathbb{G}_S}
\newcommand {\dis} {\bb{Q}_p/\bb{Z}_p}
\newcommand {\dip} {\bb{Z}/p\bb{Z}}

\newcommand {\bi} {\text{Biext}}

\title {Metric groups attached to skew-symmetric biextensions}
\author {Swarnendu Datta}
\date{}

\begin {document}

\maketitle

\begin {abstract}
Let $G$ be a commutative, unipotent, perfect, connected group scheme over an algebraically closed field of characteristic $p > 0$ and let $E$ be a biextension of $G \times G$ by the discrete group $\dis$. When $E$ is skew-symmetric, V. Drinfeld defined a certain metric group $A$ associated to $E$ (when $G$ is the perfectization of the additive group $\mathbb{G}_a$, it is easy to compute this metric group - cf. Appendix A). In this article, we prove a conjecture due to Drinfeld about the class of the metric group $A$ in the Witt group (cf. Appendix B).
\end {abstract}

\section*{Introduction}
The main result of this article (cf. thm. 1, sec. 2.3) is related to a conjecture in the theory of character sheaves on unipotent groups (developed in \cite{motivated}, \cite{foundations} as well as \cite{mitya}). Let $U$ be a unipotent group over an algebraically closed field of characteristic $p > 0$. In \cite{motivated}, the notion of $L$-packets of character sheaves on $U$ is defined in terms of certain minimal idempotents in the equivariant derived category ${\cal D}_U(U)$ (cf. \cite{motivated}, sec. 5). If $e \in {\cal D}_U(U)$ is such an idempotent, one is interested in studying the full subcategory $e{\cal D}_U(U) \subset {\cal D}_U(U)$ consisting of objects $M \in {\cal D}_U(U)$ such that $e * M \cong M$, where $*$ denotes convolution with compact support (cf. \cite{motivated}, sec. 5.5). Let ${\cal M}_e^{\text{perv}} \subset e{\cal D}_U(U)$ denote the full subcategory of $e{\cal D}_U(U)$ consisting of the perverse objects. It is shown in $\cite{foundations}$ that the following hold true (cf. Thm. 1.13):

\begin {enumerate}

\item ${\cal M}_e^{\text{perv}}$ is a semisimple abelian category with finitely many simple objects upto isomorphism (the set of isomorphism classes of simple objects is called the \emph{L-packet} associated to $e$).

\item There exists a (unique) integer $n_e$ such that $e[-n_e] \in {\cal M}_e^{\text{perv}}$. The subcategory ${\cal M}_e = {\cal M}_e^{\text{perv}}[-n_e] \subset {\cal D}_U(U)$ is closed under convolution and is a braided monoidal category with unit object $e$. Moreover, the functor $\text{id}: {\cal M}_e \rightarrow {\cal M}_e$ is equipped with a canonical automorphism $\theta$ (cf. \cite{mitya}, sec. 3.9).
\end {enumerate}
V. Drinfeld formulated the following conjecture:

\begin {conj}
$\emph{i})$ The category ${\cal M}_e$ with the structures mentioned in \emph{(2)} above is modular \emph{(}cf. \emph{sec. 2.8.2, \cite{drinfeldetel} or sec. II.1.4, \cite{turaev})}.

$\emph{ii})$ Define the functional dimension $d_e$ of ${\cal M}_e$ by  $d_e = (\dim G - n_e)/2$. The Gauss sum of the modular category \emph{(}cf. \emph{sec. 6.2, \cite{drinfeldetel})} ${\cal M}_e$ equals $\eps p^n$ where $n$ is a positive integer and $\eps = 1$ if $d_e$ is an integer and $\eps = -1$ otherwise.
\end {conj}

When $U$ is connected and $e$ is a Heisenberg idempotent in the sense of \cite{foundations} (cf. sec. 1.8), part (i) of this conjecture is proved in \cite{tanmay}. Our main theorem implies that part (ii) is also true (cf. sec. 2.4). Since the study of general minimal idempotents in ${\cal D}_U(U)$ can be reduced to the case of Heisenberg idempotents (cf. \cite{foundations}, Thm 1.2) , it is plausible that the conjecture in general can be reduced to this special case.\\

Briefly, the idea behind the proof of thm. 1 is as follows: The case of the additive group is proved by making certain explicit computations and then appealing to a descent result for skew-symmetric biextensions (cf. sec. 4). The case of arbitrary unipotent groups is then handled by reducing to the above case by means of certain reductions, explained in sec. 5. The principle behind these reductions is described in the remark at the end of sec. 6.\\

The article is organized as follows: In the first section, we briefly review Serre duality. In the second section, we review biextensions and define the metric group attached to a skew-symmetric biextension. In the third section, certain results about functoriality and descent for biextensions are proved. In the fourth section, we prove the main theorem for the case of the additive group. The fifth section is devoted to certain reductions and the proof of the theorem is completed in the sixth section. There are two appendices: The first one is devoted to explicit computations of the metric group for skew-symmetric biextensions of the additive group. The second one gives a brief introduction to metric groups.\\

{\em Acknowledgement}. I am deeply grateful to my adviser V. Drinfeld for suggesting various simplifications as well as improvements to the presentation of this article. I am indebted to D. Boyarchenko and the referee for their critical reading of the manuscript. I thank L. Breen for helpful correspondence.\\

{\em Notation}. In the remainder of this article, fix an algebraically closed field $k$ of characteristic $p > 0$. All fiber (and tensor) products are over $k$, unless stated otherwise. If $G$ is a scheme over $k$, the notation $g \in G$ means that $g$ is a $k$-point of $G$.

\section {Serre Duality}

We rapidly recall the notion of Serre duality for commutative unipotent groups, indicated in \cite{serre} and developed in \cite{begueri}.

\subsection {Definition and properties.}

Recall that a scheme $X$ in characteristic $p$ (i.e., $p$ annihilates the structure sheaf ${\cal O}_X$ of $X$) is called \emph{perfect} if the morphism ${\cal O}_X \rightarrow {\cal O}_X$ given by $f \mapsto f^p$ on the local sections of ${\cal O}_X$ is an isomorphism of sheaves.

Let $\fr{Perf}_k$ be the category of perfect schemes over $k$ and let $G$ be a commutative, unipotent, connected group scheme in $\fr{Perf}_k$. One considers the following functor, from the category $\fr{Perf}_k$ to abelian groups:
$$
S \mapsto \ext(G \times S, \dis)
$$
where $\ext(A,B)$ denotes (for commutative group schemes $A,B$) the group of isomorphism classes of extensions of $A$ by $B$:
$$
0 \rightarrow B \rightarrow C \rightarrow A \rightarrow 0
$$
where $C$ is commutative as well.

One knows that the functor defined above is representable by a group scheme $G^*$ in $\fr{Perf}_k$ (cf. \cite{begueri}, Prop. 1.2.1). The group $G^*$ is called the \emph{Serre dual} of $G$. One has the following properties (\emph{loc. cit.}):

i) $G^*$ is a connected, commutative, unipotent group scheme isogenous to $G$.

ii) The canonical homomorphism of $G$ onto $G^{**}$ is an isomorphism. 

iii) If $0 \rightarrow G^\p \rightarrow G \rightarrow G^{\p\p} \rightarrow 0$ is an exact sequence of connected, commutative, unipotent group schemes in $\fr{Perf}_k$, then so is $0 \rightarrow G^{\p\p*} \rightarrow G^* \rightarrow G^{\p*} \rightarrow 0$.\\

\emph{Remark}. The last two properties say that the (contravariant) functor $G \mapsto G^*$ is exact and involutive. Therefore, the canonical map $\Hom(G_1,G_2) \rightarrow \Hom(G_2^*, G_1^*)$ is an isomorphism.

\subsection {An example.}

Let $\padd$ be the \emph{perfectization} (cf. \cite{mitya}, appendix A.3) of the additive group $\add$. In particular, $\padd$ is the spectrum of the ring:
$$
R = k[x,x^{1/p},x^{1/p^2}, \cdots]
$$
One knows that the dual $\padd^*$ can be canonically identified with $\padd$. We recall this isomorphism: Let $S$ be an arbitrary perfect scheme and let $\adds = \padd \times S$.  The exponent of the group $\padd$ is $p$ (i.e., $p\padd = 0$), whence the group $\padd^*(S)= \ext(\adds,\dis)$ can be canonically identified with $\ext(\adds,\dip)$. If $\si$ is an endomorphism of $\adds$, let $\eps_\si$ denote the element of $\padd^*(S)$ deduced from the Lang isogeny:
$$
0 \rightarrow \dip \rightarrow \adds \xrightarrow{x \mapsto x^p-x} \adds \rightarrow 0
$$
via the morphism $\si: \adds \rightarrow \adds$. This gives a homomorphism $\en(\adds) \rightarrow \padd^*(S)$. In particular, since every element $t \in \padd(S)$ acts on $\adds$ by homothety, one obtains a homomorphism:
$$
\vp: \padd(S) \rightarrow \padd^*(S)
$$
of group functors, which one knows to be an isomorphism (cf. \cite{begueri}, Prop. 1.1.1).\\

\emph{Remark}. The above isomorphism $\vp$ gives rise to a \lq\lq universal\rq\rq\ scheme $E$ over $\padd \times \padd$ which satisfies the following property: For all $t \in \padd$, the fiber $E_t$ over $t \times \padd = \padd$ is the extension of $\padd$ given by $\vp(t)$. The scheme $E$ can be explicitly described as follows: Let $ T = R \otimes R$ and let $u,v \in T$ be the elements $x \otimes 1$ and $1 \otimes x$. Then $E = \spec(T[z]/(z^p - z -uv))$, equipped with its natural projection onto $\spec(T) = \padd \times \padd$. Note that the symmetry of the expression $z^p - z -uv$ as a function of $u,v$ implies that $\vp$ is self dual, i.e., $\vp = \vp^*$.

The scheme $E$ can also be viewed as a biextension of $\padd \times \padd$. The concept of a biextension is reviewed in the following section.

\section {Biextensions}

Biextensions were introduced by Mumford in \cite{mumford}. The definition which follows is a modification of the one that appears in \cite{mitya}. For a treatment of biextensions from the functorial point of view, we refer to Grothendieck's expos\'e VII in \cite{sga}. Note however that the definition of biextension which follows is weaker than the ones considered in \cite{mumford} or \cite{sga}, since we don't assume the partial groups laws to be commutative. Such biextensions were considered by Breen in \cite{breen2} under the name of weak biextensions.

Let $G_1$, $G_2$ be group schemes over $k$, and let $A$ be a commutative group scheme over $k$. A \emph{bi-extension} of $(G_1,G_2)$ by $A$ is a scheme $E$ equipped with the following additional structures:

\begin{enumerate}
\item[i)] $E$ is a group scheme over $G_2$ represented as a central extension of $G_1\times G_2$ by $A\times G_2$;
\item[ii)]  $E$ is a group scheme over $G_1$ represented as a central extension of $G_1\times G_2$ by $G_1\times A$.
\end{enumerate}
One requires that

\begin{enumerate}
\item[a)]  data (i) and (ii) give rise to the same morphism $\pi:E \rightarrow G_1\times G_2$ and the same action of $A$ on $E$;
\item[b)]  the group operations  $\bullet_1:E\times_{G_2}E \rightarrow E$ and $\bullet_2:E\times_{G_1}E \rightarrow E$ corresponding to (i) and (ii) satisfy the identity

\begin{equation}
\label {lot_of_bullets}
(e_{11}\bullet_2 e_{12}) \bullet_1 (e_{21}\bullet_2 e_{22}) =
(e_{11}\bullet_1 e_{21}) \bullet_2 (e_{12}\bullet_1 e_{22});
\end{equation}

More precisely, \eqref{lot_of_bullets} should hold for any $k$-scheme $S$ and any $S$-points $e_{11},e_{12},e_{21},e_{22} \in E(S)=\Hom_{k-\text{sch}}(S,E)$ whenever both sides of \eqref{lot_of_bullets} are defined, i.e., whenever
\[
\pi(e_{11})=(g_1,g_2), \quad \pi(e_{12})=(g_1,g_2'), \quad
\pi(e_{21})=(g_1',g_2), \quad \pi(e_{22})=(g_1',g_2')
\]
for some $g_1,g_1'\in G_1(S)$ and $g_2,g_2'\in G_2(S)$.
\end{enumerate}

\emph {Remark}. Either (i) or (ii) defines on $E$  a structure of $A$-torsor over $G_1\times G_2$. By~(a), these two structures are the same.\\

In the rest of this article, we shall always assume that $A$ is the discrete group $\bb{Q}_p/\bb{Z}_p$.

\subsection {Biadditive pairings attached to biextensions.}

Let $G_1, G_2$ be commutative, unipotent, connected group schemes in $\fr{Perf}_k$. Let $\bi(G_1,G_2)$ denote the group of isomorphism classes of biextensions of $G_1 \times G_2$ by the discrete group $\bb{Q}_p/\bb{Z}_p$. One has isomorphisms of abelian groups (cf. Lemma A.17, \cite{mitya}):
$$
\Hom(G_1,G_2^*) \leftrightarrow \bi(G_1,G_2) \leftrightarrow \Hom(G_2,G_1^*)
$$
where the composite map (in both directions) is the one obtained from duality. We refer to \cite{sga}, expos\'e VII.1.6 for the definition of these homomorphisms (which extend easily to the more general case considered here).

Let $E$ be a biextension of $G_1 \times G_2$ and let $f: G_1 \rightarrow G_2^*$ be the corresponding homomorphism. Then the restriction $E|_{\ker f \times G_2}$ (resp. $E|_{G_1\times \ker f^*}$) is the trivial biextension of $\ker f\times G_2$ (resp. $G_1 \times \ker f^*$). Since $G_1,G_2$ are connected, we get \emph{unique} trivializations (cf. \cite{mitya}, Cor.  A.13):
\begin {align*}
\al_1& : (\ker f \times G_2) \times \dis \xrightarrow{\sim} E|_{\ker f \times G_2}\\
\al_2& : (G_1 \times \ker f^*) \times  \dis \xrightarrow{\sim} E|_{G_1\times_k \ker f^*}
\end {align*}
Both $\al_1$ and $\al_2$ induce trivializations of $E|_{\ker f \times \ker f^*}$ by restriction, whence the composition $\al_1^{-1} \circ \al_2$ (suitably restricted) gives an automorphism of the trivial biextension of $\ker f \times \ker f^*$. Therefore we get a bi-homomorphism (of group schemes):
$$
B_0: \ker f \times \ker f^* \rightarrow \dis
$$
Since $B_0$ is trivial on $(\ker f)^0 \times (\ker f^*)^0$, one obtains a biadditive pairing
$$
B: \pi_0(\ker f) \times \pi_0(\ker f^*) \rightarrow \bb{Q}_p/\bb{Z}_p
$$
by passing to the quotient. For explicit computation of $B$ in the case when $G_1 = \padd = G_2$, see prop.~\ref{formulabiadditive}, Appendix A, as well as the remark following it.

The pairing $B$ can be made explicit as follows: Let $a \in \ker f$ and $b \in \ker f^*$. Since $G_2$ can be identified with the dual of $G_2^*$, the element $b \in G_2$ corresponds to an  extension of $G_2^*$:

$$
0 \rightarrow \dis \rightarrow \wt{G_2^*} \rightarrow G_2^* \rightarrow 0
$$
The hypothesis $b \in \ker f^*$ implies that the pullback of the above extension to $G_1$ via the homomorphism $f$ is the trivial extension. Therefore one has a (unique) factorization of $f$ via $\wt{G_2^*}$:

$$
\xymatrix{
   & & & G_1 \ar[d]^-{f} \ar@{-->}[dl]^-{\wt{f}} &\\
 0 \ar[r] & \dis \ar[r] & \wt{G_2^*} \ar[r] & G_2^* \ar[r] & 0
}
$$
As $a \in \ker f$, it follows that $\wt{f}(a)$ is an element of $\dis$. It can be checked that $B_0(a,b) = \wt{f}(a)$.\\

\emph {Remark}. If $f$ is an isogeny then one has an exact sequence:
$$
0 \rightarrow \Hom(\ker f, \dis) \rightarrow \ext(G_2^*, \dis) \rightarrow \ext(G_1, \dis)
$$
which shows that $\ker f^*$ can be canonically identified with $\Hom(\ker f, \dis)$, i.e., the Pontryagin dual of $\ker f$. It is easy to check that the pairing $B$ comes from this duality.

\begin {prop}
\label {nondegenerate}
The pairing $B$ is nondegenerate.
\end {prop}

\begin {proof}
See \cite{mitya}, Prop. A. 19.\qed
\end {proof}

\subsection {Skew-symmetric biextensions.}

We now specialize to the following situation: $G_1 = G = G_2$.

\begin {define}
\emph{i)} A biextension $E$ of $G \times G$ is \emph{skew-symmetric} if the pullback of $E$ via the diagonal homomorphism $G \xrightarrow{\dde} G\times G$ is trivial $($as a $\dis$ torsor$)$. 

\emph{ii)} A homomorphism $f: G \rightarrow G^*$ is \emph{skew-symmetric} if the corresponding biextension is skew-symmetric.
\end {define}

\emph{Remarks}. a) Since the natural map $\bi(G,G) \rightarrow \Hom(G,G^*)$ is an isomorphism of \emph{abelian groups}, it follows that skew-symmetric homomorphisms form a subgroup of $\Hom(G,G^*)$.

b) In \cite{breen}, L. Breen introduced a more general notion of a skew-symmetric biextension (called biextensions altern\'ees, cf. sec. 1, {\em loc. cit.}) by a sheaf of abelian groups.\\

Let $f: G \rightarrow G^*$ be the homomorphims corresponding to a biextension $E$ of $G \times G$. One considers the following condition:
\begin {equation}
f + f^* = 0
\tag{$*$}
\end {equation}

\begin {prop}
\label {skew}
\emph{i)} If $E$ is skew-symmetric then $(*)$ above holds. Conversely, if $p \neq 2$, then $(*)$ is equivalent to the skew-symmetry of $E$.

\emph{ii)}  If $E$ is defined by a homomorphism $G \rightarrow G^*$ of the form $w - w^*$, then $E$ is skew-symmetric.
\end {prop}

\emph{Remark}. In view of (i) above, one says that a biextension $E$ is \emph{weakly skew-symmetric} if the homomorphism $f$ satisfies the condition $(*)$. \\

\begin {proof}
The proof is standard.\qed
\end {proof}

\begin {prop}
If $E$ is skew-symmetric, then $\ker f = \ker f^*$. Further, the pairing $B$ \emph{(cf. sec. 2.1)} is \emph{symmetric}, i.e., $B(a,b) = B(b,a)$ for all $a,b \in \ker f$.
\end {prop}

\begin {proof}
The first statement follows from part (i) of the previous proposition. For the proof of the second one, see \cite{mitya}, Lemma A.25.\qed
\end {proof}

\subsection {Metric groups attached to skew-symmetric biextensions.}

Let $E$ be a skew-symmetric biextension of $G \times G$ and let $f: G \rightarrow G^*$ be the corresponding homomorphism. Set $A = \pi_0(\ker f)$, which is a finite $p$-group. We define a map
$$
q: A \rightarrow \dis
$$
as follows: By hypothesis, the $\dis$ torsor $E|_{\dde(G)}$ over $\dde(G)$ is trivial, whence there exists a \emph{unique} trivialization:
$$
\vp: G \times \dis \rightarrow E|_{\dde(G)}
$$
such that $\vp(0,0) = 0$. On the other hand, one has a trivialization:
$$
\vp^\p: \ker f \times \dis \rightarrow E|_{\dde(\ker f)}
$$
of $E|_{\dde(\ker f)}$ deduced from $\al_1$ by restriction (cf. sec. 2.1) via the embedding $\ker f \times \dis \hookrightarrow (\ker f \times G_2) \times \dis$ (where $\ker f$ is embedded in $\ker f \times G_2$ via the diagonal map). Therefore the composition $\vp^{-1} \circ \vp^\p$ (suitably restricted) gives an automorphism of the trivial $\dis$ torsor over $\ker f$, which defines a morphism (of schemes):
$$
q_0: \ker f \rightarrow \dis
$$
One obtains $q$ by passing to the quotient. For explicit computation of $q$ in the case when $G = \padd$, see prop.~\ref{formulaquadratic} in Appendix A.

\begin {prop}
\label {polarization}
The function $q$ is a quadratic form $($cf. Appendix B$)$ on $A$. Further, if $B$ is the biadditive pairing on $A \times A$ \emph{(cf. sec. 2.1)}, then
$$
B(a,b) = q(a+b)-q(a)-q(b)
$$
i.e., $B$ is the \emph{polarization} of $q$. 
\end {prop}

\begin {proof}
See \cite{mitya}, Lemma A.26.\qed
\end {proof}

In particular, the quadratic form $q$ is nondegenerate, whence the group $A$ equipped with $q$ is a \emph{metric group} (cf. Appendix B). The conjecture which appears in \cite{drinfeld} is the following:

\begin {thm} Let $G$ be a perfect, commutative, unipotent, connected group scheme over $k$ and let $E$ be a skew-symmetric biextension of $G$. Then with the notation above, one has:

$\emph{a})$ $\log_p(\emph{\text{Card}}(A))$ is even.

$\emph{b})$ Let $[(A,q)]$ denote the class of $(A,q)$ in the Witt group $W$ $($cf. appendix$)$. Then
$$
[(A,q)] =
\begin {cases}
0 & \text{if $\dim f(G)$ is even} \\
(\bb{F}_{p^2}, i\circ\text{\emph{Nm}}_{\bb{F}_{p^2}/\bb{F}_p}) & \text{if $\dim f(G)$ is odd}
\end {cases}
$$
where $i: \bb{F}_p \hookrightarrow \dis$ is the embedding defined by $i(1) = 1/p$.
\end {thm}

\emph{Remark}. Part (b) above can be reformulated as follows: If $\dim f(G)$ is even, then $(A,q)$ is \emph{hyperbolic}, i.e., admits a Lagrangian subgroup (cf. Appendix B). Otherwise, there exists an isotropic subgroup (cf. Appendix B) $I \subset A$ such that $I^\perp/I \cong (\bb{F}_{p^2}, i\circ \text{Nm}_{\bb{F}_{p^2}/\bb{F}_p})$.\\

The theorem is proved for $G = \padd$ in sec. 4. The general case, which is treated in the last two sections, is proved by reducing to the above case.

\subsection {An application}

In this section, we prove part (ii) of the conjecture mentioned in the introduction when $U$ is a connected unipotent group and $e$ is a Heisenberg idempotent. We recall that such idempotents are defined in terms of (cf. \cite{mitya}, sec. 6.1):

a) A connected normal subgroup $N$ of $U$ such that $U/N$ is commutative.

b) A $U$-equivariant multiplicative local system on $N$ (here $U$ acts on $N$ by conjugation).

Such a local system gives rise to a skew-symmetric biextension of $(U/N)_{\text{perf}}$, i.e., the perfectization of $U/N$. One requires that the corresponding homomorphism $(U/N)_{\text{perf}} \rightarrow (U/N)_{\text{perf}}^*$ be an isogeny.

Let $A$ be the metric group corresponding to this biextension. The following result has been proved in \cite{tanmay}:

\begin {thm}
The category ${\cal M}_e$ is equivalent to the modular category corresponding to the metric group $A$.
\end {thm}

Since the Gauss sum of a modular category (cf. \cite{drinfeldetel}, sec. 6) corresponding to a metric group is none other than the Gauss sum of the associated metric group, it follows that part (ii) of the conjecture can be rephrased in terms of the Gauss sum of $A$. Note also that the integer $n_e$ equals $\dim (N)$ (cf. \cite{foundations}, Thm. 1.26/(b)) so that the functional dimension $d_e$ is $1/2 \dim(U/N)$. Therefore part (ii) immediately follows from the following (by taking $G = (U/N)_{\text{perf}}$):

\begin {prop}
Let $G$ be a commutative, unipotent, perfect, connected group scheme over $k$ and $f: G \rightarrow G^*$ be a skew-symmetric homomorphism, which is also an isogeny. Then the Gauss sum of the metric group $A = \ker f$ equals $\eps p^n$ where $n$ is an integer and $\eps = 1$ if $\dim G$ is even and $-1$ otherwise.
\end {prop}

\begin {proof}
We shall need the following fact (cf. \cite{drinfeldetel}, Prop. 6.1): Let $A$ be a metric group and $H \subset A$ be an isotropic subgroup (cf. Appendix B). Let $B$ be the group $H^\perp/H$ with the induced metric group structure. If $\tau(A)$ (resp. $\tau(B)$) denotes the Gauss sum of $A$ (resp. $B$), then
$$
\tau(A) = \card(H) \tau(B)
$$
Therefore the Gauss sums of two metric groups which are Witt equivalent (cf. Appendix B) differ by a power of $p$. By virtue of thm. 1, it remains to note that:

a) $\tau(A) = 1$ when $A = 0$.

b) $\tau(A) = -p$, when $A = (\bb{F}_{p^2}, i\circ\text{Nm}_{\bb{F}_{p^2}/\bb{F}_p})$.

\noindent The verification of (b) is easy and is carried out in \cite{drinfeldetel}, Prop. A.3/(ii).\qed
\end {proof}

\section {General Results}
In this section we state several results regarding functoriality (sec. 3.1) and descent (sec. 3.2) of biextensions. In the following, $G_1, G_2$ will denote commutative, unipotent, perfect, connected groups.

\subsection {Functoriality.}
We explicitly state two (obvious) functorial results about the pairing $B$ and the form $q$ defined in sec. 2.1 and sec. 2.3.

Let $E$ be a biextension of $G_1 \times G_2$ corresponding to $f: G_1 \rightarrow G_2^*$. Let $\vp_1: G_1^\p \rightarrow G_1$ and $\vp_2: G_2^\p \rightarrow G_2$ be homomorphisms and $E^\p$ be the biextension of $G_1^\p \times G_2^\p$ obtained from $E$ by pulling back via $\vp_1 \times \vp_2$. 

Let $a^\p \in G_1^\p$, $b^\p \in G_2^\p$ and $a,b$ be their images in $G_1, G_2$. Suppose that $a \in \ker f$ and $b \in \ker f^*$. This implies $a^\p \in \ker f^\p$, $b^\p \in \ker f^{\p*}$.

\begin {prop}
\label {functoriality}
Under the above assumptions, one has $B^\p(a^\p,b^\p) = B(a,b)$ where $B^\p, B$ are the pairings attached to $E^\p, E$ resp.
\end {prop}

\emph{Remark}. We are identifying $B^\p,B$ with $B^\p_0,B_0$, cf. the definition of $B$.

\begin {proof}
This follows immediately from the definition of the pairings $B^\p, B$.\qed
\end {proof}
Let $\vp: G^\p \rightarrow G$ be a homomorphism (with $G^\p$ connected) and let $E^\p$ be the biextension of $G^\p$ deduced from a skew-symmetric biextension $E$ via $\vp$. Let $a^\p \in G^\p$ such that $\vp(a^\p) = a \in \ker f$. Then (with obvious notation):

\begin {prop}
\label{quadfunctoriality}
One has $q^\p(a^\p) = q(a)$.
\end {prop}

\begin {proof}
This follows immediately from the definition of $q^\p, q$.\qed
\end {proof}

\emph{Remark}. Suppose $G^\p \xrightarrow{\vp} G$ is surjective with \emph{connected} kernel $H$. Note that $E^\p$ corresponds to the homomorphism $f^\p: G^\p \rightarrow G^{\p*}$ which is the composition:
$$
G^\p \xrightarrow{\vp} G \xrightarrow{f} G^* \xrightarrow{\vp^*} G^{\p*}
$$
Since $\vp^*$ is injective, one has $\vp^{-1}(\ker f) = \ker f^\p$, i.e., $\vp(a^\p) \in \ker f$ if and only if $a^\p \in \ker f^\p$.\\

\subsection {Descent.}

We keep the notation of the last section. Let $H_1$ (resp. $H_2$) be a closed subgroup of $G_1$ (resp. $G_2$). We do not assume $H_1, H_2$ to be connected.

\begin {prop}
\label {descent}
The following conditions are equivalent:

\emph{i}$)$ $E$ is obtained from a biextension of $G_1/H_1 \times G_2/H_2$ by pullback via the natural projection $G_1 \times G_2 \surject G_1/H_1 \times G_2/H_2$.

\emph{ii}$)$ The following are satisfied:

\emph{a}$)$ $H_1 \subset \ker f$ and $H_2 \subset \ker f^*$.

\emph{b}$)$ $B(H_1,H_2) = 0$, i.e., $B(a,b) = 0$ for all $a \in H_1, b \in H_2$.
\end {prop}

\begin {proof}

(i) $\implies$ (ii): Condition (a) is clear and (b) follows from prop.~\ref{functoriality}.

(ii) $\implies$ (i): Since $H_1 \subset \ker f$, one has a section $s_1: H_1 \times G_2 \rightarrow E$. This section allows one to view $H_1 \times G_2$ as a group subscheme of $E$, where $E$ is considered as a group scheme over $G_2$. Thus the group $H_1 \subset H_1 \times G_2$ acts by translation on $E$, preserving the projection $E \rightarrow G_2$. Similarly, as $H_2 \subset \ker f^*$, one has a section $s_2: G_1 \times H_2 \rightarrow E$ which embeds $G_1 \times H_2$ as a group subscheme of $E$, where $E$ is considered as group scheme over $G_1$. As above, one obtains an action of $H_2$ on $E$. One has:

\begin {enumerate}

\item In virtue of (ii)/(b), the sections $s_1$ and $s_2$ agree on $H_1 \times H_2$. This implies that the actions of $H_1$ and $H_2$ on $E$ commute with each other (this follows from the compatibility condition in the definition of biextensions). Thus one obtains an action of $H_1 \times H_2$ on $E$.

\item The projection $E \rightarrow G_1 \times G_2$ is equivariant under the $H_1 \times H_2$ action, where $H_1 \times H_2$ acts on $G_1 \times G_2$ by translation.

\end {enumerate}
Let $E^\p = E/(H_1 \times H_2)$. In view of above, one obtains a morphism $\pi^\p: E^\p \rightarrow G_1/H_1 \times G_2/H_2$, which gives $E^\p$ a structure of a biextension over $G_1/H_1 \times G_2/H_2$ (the group scheme structures of $E^\p$ over $G_1/H_1$ and $G_2/H_2$ are deduced from those of $E$ over $G_1$ and $G_2$ by factoring out by $H_1$ and $H_2$-actions) , from which $E$ is deduced via pullback, qed.\qed
\end {proof}

\emph{Remarks}. 1) When $H_2$ is connected, the group $\Hom(H_2,\dis)$ is reduced to zero, whence condition (b) in (ii) is vacuous. The same remark applies when $H_1$ is connected. 

2) In particular, if $H_1 = \ker f$ and $H_2 = \ker f^*$, then the biextension $E$ of $G_1 \times G_2$  is deduced from a biextension of $G_1/H_1^0\times G_2/H_2^0$. Note that the corresponding homomorphism $G_1/H_1^0 \rightarrow (G_2/H_2^0)^*$ is an isogeny. Indeed one observes that this map as well as its dual has finite kernel and uses property (i) and (iii), sec. 1.1.

\begin {prop}
\label {descentskewsymmetric}
Let $E$ be a skew-symmetric biextension of $G \times G$ and let $A_0 = \ker f$. If $H$ is a closed subgroup of $A_0$ such that $q_0|_{H} = 0$ \emph{(cf. sec. 2.3)}, then $E$ is deduced from a skew-symmetric biextension of $G/H \times G/H$.
\end {prop}

\begin {proof}
Since $H \subset \ker f$ and $B(H,H) = 0$ (in virtue of prop.~\ref{polarization}), it follows from prop.~\ref{descent} that $E$ is deduced from a biextension $E^\p$ of $G/H \times G/H$. It suffices to show that $E^\p$ is skew-symmetric. Since $E$ is skew-symmetric, there exists a (unique) morphism $G \xrightarrow{u} E$ which takes 0 to 0 and such that the composite $G \xrightarrow {u} E \rightarrow G \times G$ is the diagonal embedding. We want to descent $u$ to a morphism $G/H \xrightarrow{u^\p} E^\p = E/(H \times H)$. For this, it suffices to show that the composite $G \xrightarrow{u} E \rightarrow E^\p$ is $H$-equivariant, where $H$ acts on $G$ by translation and on $E^\p$ trivially. Note that if we view $H$ as a subgroup scheme of $H \times H$ via the diagonal embedding, we obtain an action of $H$ of $E$. The projection $E \rightarrow E^\p$ is $H$ equivariant, whence it suffices to show that $G \xrightarrow{u} E$ is also $H$ equivariant. This follows from the hypothesis that $q_0|_{H} = 0$.\qed
\end {proof}

{\em Remark}. In \cite{breen}, Prop. 2.3, Breen formulates a descent result for skew-symmetric biextensions in a more general situation.

\section {The case of $\padd$}
We shall use the following notation (cf sec. 1.2):
\begin {align*}
R & = k[x,x^{1/p}, x^{1/p^2},\cdots]\\
T & = R \otimes R\\
u &= x \otimes 1 \in T,\quad v = 1 \otimes x \in T
\end {align*}
as well as 
$$
\vp: \padd \rightarrow \padd^*
$$
to denote the identification of $\padd$ with its dual $\padd^*$, as given in sec. 1.2. 

\subsection {Skew-symmetric homomorphisms.}
We want to give an explicit description of all skew-symmetric homomorphisms from $\padd$ onto its dual. All endomorphisms (of $\padd$, $\padd^*$ etc.) considered below are endomorphisms in the sense of group schemes.

If $\si: \padd \rightarrow \padd$ is an endomorphism, then one has a corresponding endomorphism $\si^*: \padd^* \rightarrow \padd^*$. The identification of $\padd$ with $\padd^*$ via $\vp$ therefore gives an anti-automorphism:
$$
\en(\padd) \rightarrow \en(\padd)
$$
This can be made explicit as follows:

\begin {lem}
\label{skewhom}
\emph{i)} $\emph{End}(\padd) = k\{\tau,\tau^{-1}\}$ where $\tau$ is the Frobenius automorphism $($which sends $x$ to $x^p$$)$. The algebra structure is given by $\tau a = a^p\tau$ for all $a \in  k$.

\emph{ii)} The map $\si \mapsto \si^*$ is given by $\tau^* = \tau^{-1},\ c^* = c$, for all $c \in k$.
\end {lem}

\begin {proof}

i) If $\si$ is an endomorphism of $\padd$, one can identify $\si$ with an endomorphism of the ring $R$. Then $\si(x)$ is a polynomial in $x^{1/p^n}$ for sufficiently large $n$. It follows that $\si\tau^n(x)$ is a polynomial in $x$. One easily checks that this polynomial is additive, whence the endomorphism $\si\tau^n$ of $\padd$ can be written as a polynomial in $\tau$ (with coefficients in $k$). This proves that $\en(\padd) = k\{\tau,\tau^{-1}\}$. The commutation relation $\tau a = a^p\tau$ is easy to verify.

ii) Let $\eps$ be the Lang isogeny of $\padd$. If $\si \in \en(\padd)$, let $\si^+(\eps)$ be the element of $\padd^*$ deduced from $\eps$ by pulling back via the morphism $\si$. If we identify $b \in k$ with the corresponding homothety of $\padd$, then recall that the element $b^+(\eps) \in \padd^*$ corresponds to the element $b \in \padd$ via $\vp$. The equality $\tau^* = \tau^{-1}$ thus amounts to showing that the pullback of the element $a^+(\eps)$ of $\padd^*$ via $\tau$ is isomorphic to $(\tau^{-1}(a))^+(\eps) = (a^{1/p})^+(\eps)$. But
$$
\tau^+a^+(\eps) = (a\tau)^+(\eps) = (\tau a^{1/p})^+(\eps) = (a^{1/p})^+\tau^+(\eps)
$$
whence it remains to check that $\tau^+(\eps) = \eps$. Let
$$
0 \rightarrow \dip \rightarrow \wt{\padd} \rightarrow \padd \rightarrow 0
$$
be the extension corresponding to $\tau^+(\eps)$. One checks that the morphism $\wt{\padd} \rightarrow \padd$ is obtained from the ring homomorphism:
$$
R \rightarrow R[z]/(z^p -z - x^p)
$$
Setting $z^\p = z-x$, the above inclusion can be written as $R \rightarrow R[z^\p]/(z^{\p p}-z^\p-x)$, which corresponds to $\eps$, qed.\qed
\end {proof}

Let $V$ be the subgroup of $k\{\tau,\tau^{-1}\}$ generated by elements of the form $c\tau^j-\tau^{-j}c$. 

\begin {cor}
\label {weakskewforga}
A homomorphism $f: \padd \rightarrow \padd^*$ is weakly skew-symmetric $($i.e., $f+f^* = 0)$ if and only if:

\emph {i)} $f \in V$, provided $p \neq 2$.

\emph {ii)} $f \in V \oplus k$ if $p=2$.
\end {cor}

\emph {Remark}. We are identifying the sets $\en(\padd)$ and $\Hom(\padd,\padd^*)$ by sending $f \in \en(\padd)$ to $\vp \circ f$.

\begin {proof}
Indeed, this is clear in view of part (ii) of the lemma above.\qed
\end {proof}

\begin {cor}
\label {skewforga}
A homomorphism $f: \padd \rightarrow \padd^*$ is skew-symmetric if and only if $f \in V$.
\end {cor}

\begin {proof}
When $p$ is odd, this follows from prop.~\ref{skew}/(i), sec. 2.2 and the above cor.~\ref{weakskewforga}. Assume then $p = 2$. Since skew-symmetric homomorphisms form a subgroup of $k\{\tau,\tau^{-1}\}$ (cf. remark following the definition of skew-symmetry), it suffices to show the following:

i) The homomorphism $c\tau^n - \tau^{-n}c$ is skew-symmetric.

ii) The homomorphism defined by a nonzero scalar $t \in k$ is not skew-symmetric. 

To see (i), apply  prop.~\ref{skew}/(ii) to $w = c\tau^n$. For (ii), denote by $E$ the corresponding biextension. One checks easily that
$$
E = \spec(T[z]/(z^2-z-tuv))
$$
Therefore if $E_\dde$ is the $\bb{Z}/2\bb{Z}$ torsor over $\padd$ obtained by pulling back $E$ via the diagonal map, one has
$$
E_\dde = \spec(R[z]/(z^2-z-tx^2))
$$
It follows that $E_\dde$ corresponds to the nonzero element $t^{1/2}$ of $\padd^*$ and thus nontrivial as a torsor, qed.\qed
\end {proof}

\subsection {The proof of Theorem 1 for $G = \padd$}

Fix a skew-symmetric homomorphism $f: \padd \rightarrow \padd^*$. If $f = 0$, then $A = \pi_0(\ker f) = 0$, whence the result becomes trivial. Suppose then $f \neq 0$. Then cor.~\ref{skewforga} allows one to write
$$
f = c_n\tau^n - \tau^{-n}c_n + \sum_{1 \leq j < n}(c_j\tau^j - \tau^{-j}c_j)
$$
with $c_n \neq 0$. The group $\ker f$ is finite with cardinality $p^{2n}$, whence part (a) follows (since $A = \ker f$). 

To prove (b), we shall use the following (easy) fact:

\begin {lem}
\label{metrickilledbyp}
A metric group of order $p^{2n}, n \in \bb{N}$ whose associated quadratic form takes values in $p^{-1}\bb{Z}_p/\bb{Z}_p \subset \dis$ is either hyperbolic or Witt equivalent to $(\bb{F}_{p^2}, i\circ\text{\emph{Nm}}_{\bb{F}_{p^2}/\bb{F}_p})$ where $i: \bb{F}_p \hookrightarrow \dis$ is defined by $i(1) = 1/p$.
\end {lem}

\begin {proof}
This follows from \cite{drinfeldetel}, Prop. A.3 and the observation that every metric group is Witt equivalent to an anisotropic one (deduced by passing to the subquotient corresponding to a maximal isotropic subgroup, cf. Appendix B).\qed
\end {proof}

Note that since $\padd$ has exponent $p$, the metric group $(A,q)$ takes values in $\bb{F}_p$ = $p^{-1}\bb{Z}_p/\bb{Z}_p$, whence in virtue of lem.~\ref{metrickilledbyp}, it suffices to prove that $(A,q)$ is not hyperbolic, i.e., does not admit a Lagrangian subgroup (cf. Appendix B). Assume to the contrary that $(A,q)$ is hyperbolic. Let $L$ be an isotropic (cf. Appendix B) subgroup of $(A,q)$. It follows from prop.~\ref{descentskewsymmetric}, sec. 3.2, that the biextension $E$ corresponding to $f$ is deduced from a skew-symmetric biextension $E^\p$ of $\padd/L \cong \padd$. Let $f^\p: \padd \rightarrow \padd^*$ be the homomorphism corresponding to $E^\p$ and let $A^\p = \card(\ker f^\p)$. We claim that:
$$
\card(A) = \card(A^\p)\card(L)^2
$$
Indeed, since $f$ is the composite $\padd \xrightarrow{\pi} \padd/L \xrightarrow{f^\p} (\padd/L)^* \xrightarrow{\pi^*} \padd^*$, it suffices to prove that $\card(L) = \card(\ker \pi^*)$. This follows from:
\begin {lem}
\label {equal}
If $u: G^\p \rightarrow G$ is an isogeny, then $\emph{Card}(\ker u) = \emph{Card}(\ker u^*)$.
\end {lem}
\begin {proof}
This is a consequence of the nondegeneracy of the biadditive pairing associated to $u$, cf. prop.~\ref{nondegenerate}, sec. 2.1. Alternatively, one could use the exact sequence in the remark preceding prop.~\ref{nondegenerate} to directly deduce that the groups $\ker u$ and $\ker u^*$ are dual to each other.\qed
\end {proof}
Assume now that $L$ is \emph{Lagrangian} (cf. Appendix B). It follows that $\card(A^\p) = 1$, i.e., $f^\p$ is an isomorphism. Since $f^\p$ is skew-symmetric, cor.~\ref{skewforga} gives a contradiction, which completes the proof.

\section {The General Case: Reductions}

Let now $G$ be a commutative, unipotent, perfect, connected group and $f: G \rightarrow G^*$ be a skew-symmetric homomorphism corresponding to a biextension $E$ of $G \times G$. One easily reduces to the case when $f$ is an isogeny. Indeed, in any case, $E$ induces a skew-symmetric biextension of $G/(\ker f)^0 \times G/(\ker f)^0$ (cf. prop~\ref{descentskewsymmetric}, sec. 3.2) and the corresponding homomorphism of $G/(\ker f)^0$ onto its dual is an isogeny (cf. remark 2 following prop.~\ref{descent}, sec. 3.2). It now suffices to use the following (with $G^\p = G/(\ker f)^0$), where recall that $A = \pi_0(\ker f)$ is equipped with the quadratic form $q$ (cf. sec. 2.3):

\begin {lem}
\label{skewdescent}
Let $\pi:G \rightarrow G^\p$ be a surjective homomorphism with \emph{connected} kernel. If $E$ descends to a biextension $E^\p$ on $G^\p$ $($i.e., $E$ is deduced from $E^\p$ by pulling back via $\pi \times \pi)$, then:

\emph{i)} $E^\p$ is a skew-symmetric biextension.

\emph{ii)} If $(A^\p,q^\p)$ is the metric group attached to $E^\p$, one has a canonical isomorphism $(A,q) \cong (A^\p,q^\p)$.
\end {lem}

\begin {proof}
i) This follows immediately from prop.~\ref{descentskewsymmetric}, sec. 3.2.

ii) It is clear that $A$ and $A^\p$ are canonically identified. To see that this identification gives a morphism of metric groups, one applies prop.~\ref{quadfunctoriality}, sec. 3.1, as well as the remark following it.\qed
\end {proof}

We shall assume $f$ is an isogeny (and thus $A = \ker f$) in what follows.

Let $H$ be a closed connected subgroup of $G$. We define $H^\perp$ to be the kernel of the following composition:
$$
G \xrightarrow{f} G^* \surject H^*
$$
One says that $H$ is \emph{isotropic} if $H \subset H^\perp$. On the other hand, if $H \cap H^
\perp$ is finite, one says that $H$ is \emph{nondegenerate}.

\subsection {The isotropic case.} 

Let $H$ be an isotropic subgroup of $G$ and let $H^\p = (H^\perp)^0$. Denote by $f^\p$ the following composition:
$$
H^\p \hookrightarrow G \xrightarrow{f} G^* \surject H^{\p*}
$$
Note that $f^\p$ corresponds to the biextension of $H^\p$ deduced from $E$ by restriction. 

\begin {lem}
\label {component}
The neutral connected component of $\ker f^\p$ is $H$.
\end {lem}

\begin {proof}
Indeed, it is clear that $H \subset (\ker f^\p)^0$ since $H$ is isotropic. For the opposite inclusion, it suffices to show both groups have same dimension. For this, note that the map $(\ker f^\p)^0 \rightarrow G^*$ (obtained by restricting $f$) admits a factorization $(\ker f^\p)^0 \xrightarrow{g} (G/H^\p)^* \hookrightarrow G^*$. Since $\ker g \subset \ker f$, it follows that $\dim (\ker f^\p)^0 \leq \dim (G/H^\p)^* = \dim G/H^\p = \dim H$, whence the result.\qed
\end {proof}

We introduce the following notation:
\begin {align*}
A_H & = A \cap H\\
A_{H^\p} & = A \cap H^\p \\
A_{H^\p/H} & = (\ker f^\p)/(\ker f^\p)^0\\
q^\p & = \text{ quadratic form on $A_{H^\p/H}$ induced by $f^\p$}
\end {align*}
Since $A = \ker f \subset \ker f^\p$, there is a natural map $A_{H^\p} \rightarrow A_{H^\p/H}$. By the above lemma, the kernel of this map is $A_H$, whence we get an injective homomorphism
$$
\vp: A_{H^\p}/A_H \rightarrow A_{H^\p/H}
$$
by passing to the quotient. One has:

i) $q(A_H) = 0$. This can be seen by applying prop.~\ref{quadfunctoriality}, sec. 3.1 to the morphism $H \hookrightarrow G$, and noting that the biextension of $H \times H$ obtained from $E$ by restriction is trivial.

ii) $A_{H^\p} \subset A^\perp_H$, i.e., if $B$ is the biadditive pairing on $A$, then $B(a,b) = 0$ whenever $a \in H^\p$ and $b \in H$. This follows from prop.~\ref{functoriality}, sec. 3.1, applied to the morphism $H^\p\times H^\p \hookrightarrow G\times G$.

Thus the group $A_{H^\p}/A_H$ is equipped with a quadratic form induced by $q$. Note that $\vp$ is a homomorphism respecting the quadratic form structures, because so is the map $A_{H^\p} \rightarrow A_{H^\p/H}$ (where $A_{H^\p}$ is equipped with $q$). Indeed, one applies prop.~\ref{quadfunctoriality}, sec. 3.1 to the morphism $H^\p \hookrightarrow G$.

\begin {lem}
\label {part1}
The homomorphism $\vp$ is an isomorphism of metric groups, where $A_{H^\p/H}$ is equipped with the form $q^\p$. 
\end {lem}

\begin {proof}
In view of the above, it suffices to show that $\vp$ is surjective. Note that the homomorphism $H \rightarrow (G/H^\p)^*$ is an isogeny (cf. the proof of lem.~\ref{component}) and in particular surjective. Thus if $a \in \ker f^\p$, there exists $b \in H$ such that $a$ and $b$ map to the same element in $(G/H^\p)^*$. It follows that $a-b$ maps to zero in both $H^{\p*}$ and $(G/H^\p)^*$, whence belongs to $\ker f$. This proves the surjectivity of $\vp$, as desired.\qed
\end {proof} 

\begin {lem}
\label {part2}
The group $A_{H^\p}/A_H$ equipped with the quadratic form induced by $q$ is a subquotient $($cf. Appendix B$)$ of $A$ $($equipped with $q)$.
\end {lem}

\begin {proof}
In other words, we want to show that the inclusion noted above
$$
A_{H^\p} \subset A_H^{\perp}
$$
is bijective. Let $a \in A$ such that $B(a,A_H) = 0$. Consider the biextension of $G \times H$ obtained from $E$ by restriction. Let $g: G \rightarrow H^*$ be the corresponding homomorphism (which is the composition $G \xrightarrow{f} G^* \rightarrow H^*)$. Then $(\ker g)/(\ker g)^0 = H^\perp/H^\p$ and $\ker g^* = A_H$. One knows that the corresponding pairing between between these groups is nondegenerate (cf. prop.~\ref{nondegenerate}, sec. 2.1). Since $a$ can be interpreted as an element of $H^\perp$, the hypothesis $B(a,A_H) = 0$ implies that $a$ is zero in $H^\perp/H^\p$, i.e., $a \in H^\p$, as desired.\qed
\end {proof}

\begin {cor}
\label {hyperbolic1}
The metric group $(A,q)$ is hyperbolic if and only if $(A_{H^\p/H}, q^\p)$ is.
\end {cor}

\begin {proof}
This follows by combining lem.~\ref{part1} and ~\ref{part2}.\qed
\end {proof}

\subsection {The nondegenerate case.}

We begin with the following general situation: Let $\vp: G^\p \rightarrow G$ be an isogeny (with $G^\p$ connected) and let $E^\p$ be the biextension of $G^{\p}$ obtained from $E$ by pulling back via the above map. Let $f^\p: G^\p \rightarrow G^{\p*}$ be the corresponding homomorphism, obtained as the composition:
$$
G^{\p} \xrightarrow{\vp} G \xrightarrow{f} G^* \xrightarrow{\vp^*} G^{\p*}
$$
Let $A^\p = \ker f^\p$ and let $q^\p$ be the quadratic form on $A^\p$. 

\begin {prop}
\label {isogeny}
$(A,q)$ is a subquotient $($cf. Appendix B$)$ of $(A^\p,q^\p)$.
\end {prop}

\begin {proof}
Note that $\vp^{-1}(A)$ is contained in $A^\p$. Thus if $K^\p = \vp^{-1}(A)$ and $K = \vp^{-1}(0)$, then one has a filtration of the group $A^\p$:
$$
0 \subset K \subset K^\p \subset A^\p
$$
such that $A = K^\p/K$. It thus suffices to prove that 

i) $q^\p(K) = 1$ and

ii) $K^\perp = K^\p$

\noindent in view of prop.~\ref{quadfunctoriality}, sec. 3.1 (applied to $\vp$). Indeed, (i) is obvious since $\vp(K) = 0$ (cf. prop.~\ref{quadfunctoriality}). Also, the inclusion $K^\p \subset K^\perp$ follows from prop.~\ref{functoriality}, sec. 3.1 (applied to $\vp \times \vp$), whence the opposite inclusion would follow if we can show that
$$
\card(K)\times \card(K^\p) = \card(A^\p)
$$
in view of the nondegeneracy of $q^\p$. The above equality can be written in the following form:
$$
\card(A^\p)/\card(A) = \card(K)^2
$$
But if $K^* = \ker(G^* \xrightarrow{\vp^*} G^{\p*})$, then the left side equals $\card(K)\times \card(K^*)$. Finally, the equality $\card(K) = \card(K^*)$ follows from lem.~\ref{equal}.\qed
\end {proof}

\emph{Remark}. The above result remains true more generally when $\vp$ is a surjection. Indeed, setting $G^{\p\p} = G^\p/(\ker \vp)^0$, the morphism $\vp$ factorizes as $G^\p \rightarrow G^{\p\p} \rightarrow G$, whence we get the result by applying lem.~\ref{skewdescent} to $G^\p \rightarrow G^{\p\p}$ and the proposition above to the isogeny $G^{\p\p} \rightarrow G$.\\

We want to apply the above result to the following situation: Let $H$ be a nondegenerate subgroup of $G$. Set $G^\p = H \times H^\p$ (recall that $H^\p = (H^\perp)^0$). Then the multiplication homomorphism
$$
G^\p \rightarrow G
$$
is an isogeny, whence $(A,q)$ is a subquotient of $(A^\p,q^\p)$.

\begin {lem}
Let $E_H$ $($resp. $E_{H^\p})$ be the biextension of $H$ $($resp. $H^\p)$ obtained from $E$ by restriction. Then $E^\p = E_H \times E_{H^\p}$.
\end {lem}

\begin {proof}
Clearly, $E_H, E_{H^\p}$ can equally be obtained from $E^\p$ by restriction. Now note that the homomorphism $G^\p \rightarrow G^{\p*}$ gives homomorphisms $H \rightarrow H^{\p*}$ and $H^\p \rightarrow H^*$. The second one is trivial by hypothesis, hence so is the first one, in view of skew-symmetry. Thus $G^\p \rightarrow G^{\p*}$ is product of the restrictions $H \rightarrow H^*$ and $H^\p \rightarrow H^{\p*}$, qed.\qed
\end {proof}

Let $(A_H,q_H)$ and $(A_{H^\p},q_{H^\p})$ be the metric groups corresponding to $E_H$ and $E_{H^\p}$.

\begin {cor}
\label {hyperbolic2}
$(A,q)$ is a subquotient $($cf. Appendix B$)$ of $(A_H,q_H) \oplus (A_{H^{\p}},q_{H^{\p}})$.
\end {cor}

\begin {proof}
Indeed, the metric group corresponding to $E^\p$ is $(A_H,q_H) \oplus (A_{H^{\p}},q_{H^{\p}})$, in view of the lemma above.\qed
\end {proof}

\section {Completion of Proof}

Let $f:G \rightarrow G^*$ be a skew-symmetric homomorphism. As explained in the beginning of the previous section, we may assume $f$ is an isogeny. We argue by induction on $n = \dim G$ to prove both (a) and (b). The case $n = 1$, i.e., $G \cong \padd$ has already been treated. Let then $n > 1$. Choose a closed subgroup $H$ of $G$ isomorphic to $\padd$ and let $g$ be the composition:
$$
H \hookrightarrow G \xrightarrow{f} G^* \surject H^*
$$
Suppose first $g = 0$. Then cor.~\ref{hyperbolic1}, sec. 5.1, reduces us to the case of $H^\p/H$, whence induction applies. If $g \neq 0$, then $g$ is an isogeny, whence we are done by  cor.~\ref{hyperbolic2}, sec. 5.2.\\

\emph{Remark}. The reduction process is similar in spirit to the following \lq\lq reduction\rq\rq:
Let $V$ be a vector space over a field equipped with a nondegenerate quadratic form $q$. If $v \in V$ is a nonzero element and $H$ is the line in $V$ spanned by $v$, then:

(i) If $q(v) \neq 0$, then one has a decomposition $V = H \oplus H^{\perp}$ of quadratic modules;

(ii) If $q(v) = 0$, then $q$ induces a nondegenerate quadratic form $q^\p$ on $H^\perp/H$. Further, the image of the forms $q$ and $q^\p$ in the Witt ring are equal.\\

\section* {Appendix A: Formulas related to biextensions of $\padd \times \padd$}

Let $E$ be a biextension of $\padd \times \padd$. In this section we indicate how to obtain explicit formulas for the pairing $B$ (cf. sec. 2.1) as well as the quadratic form $q$ (when $E$ is skew-symmetric, cf. sec. 2.3) associated to $E$.

Let $f: \padd \rightarrow \padd^*$ be the homomorphism corresponding to $E$. Using the isomorphism $\vp: \padd \rightarrow \padd^*$, we view $f$ as an element of $\en(\padd) = k\{\tau, \tau^{-1}\}$. For $t \in \padd(k)$, let $E_{t,\_}$ denote the pullback $E$ via the inclusion $\padd = t \times \padd \hookrightarrow \padd \times \padd$. Then the central extension $E_{t,\_} \in \padd^*$ of $\padd$ corresponds to the element $f(t) \in \padd(k) = k$. Thus $E_{t,\_} \cong \spec(R[z]/(z^p -z - f(t)x))$ as a scheme over $\spec(R) = \padd$. This allows one to deduce that
\begin {equation}
E \cong \spec(T[z]/(z^p - z - f(u)v))
\tag{1}
\end {equation}
as a scheme over $\spec(T) = \padd \times \padd$. Similarly, if $f^*$ is the homomorphism dual to $f$ and $E_{\_,t}$ denotes the pullback of $E$ via $\padd = \padd \times t \hookrightarrow \padd \times \padd$, then $E_{\_,t} \cong \spec(R[z]/(z^p-z-f^*(t)x))$ and 
\begin {equation}
E \cong \spec(T[z]/(z^p - z -uf^*(v)))
\tag{2}
\end {equation}

\begin {prop}
\label {formulabiadditive}
$\emph{a})$ There exists a unique element $g(u,v) \in T$ such that $g(0,0) = 0$ and
$$
g(u,v)^p - g(u,v) = f(u)v - uf^*(v)
$$

$\emph{b})$ For $a \in \ker f$, $b \in \ker f^*$, the pairing $B$ is given by the formula $B(a,b) = g(a,b)$.
\end {prop}

\begin {proof}
a) Indeed, for $h_1, h_2 \in T$, write $h_1 \sim h_2$ if $h_1 - h_2$ has the form $g^p - g$ for some $g \in T$ such that $g(0,0) = 0$. Then $au^{p^i}v^{p^j} \sim a^pu^{p^{i+1}}v^{p^{j+1}}$, from which it easily follows that $f(u)v \sim uf^*(v)$ (one uses lem.~\ref{skewhom}/(ii)). This proves the existence of $g$. The uniqueness is clear.

b) It follows that the isomorphisms (1) and (2) above are related to each other via the change of variable $z \mapsto z - g(u,v)$, whence the formula.\qed
\end {proof}

\emph{Remark}. It is easy to give an explicit formula for the element $g$ above: If $f = a\tau^n$, then one has:
\begin {equation*}
g(u,v) = 
\begin {cases}
a^{p^{-1}}u^{p^{n-1}}v^{p^{-1}} + \cdots + a^{p^{-n}}uv^{p^{-n}}  & \text{if $n > 0$}\\
0 & \text{if $n = 0$}\\
-au^{p^n}v - \cdots - a^{p^{-n-1}}u^{p^{-1}}v^{p^{-n-1}} & \text{if $n < 0$}
\end {cases}
\end {equation*}
For general $f$, one notes that $g$ varies additively with respect to $f$.\\

\emph{Example}. Suppose that $f = \tau - \tau^{-1} = -f^*$. Then $A = \ker f$ identifies with roots of $x^{p^2} - x$ in $k$, i.e., with $\bb{F}_{p^2}$. One checks that $g(u,v) = uv^{1/p} + u^{1/p}v$, whence $B: \bb{F}_{p^2} \times \bb{F}_{p^2} \rightarrow \bb{Z}/p\bb{Z}$ identifies with the pairing $(a,b) \mapsto \text{tr}(ab^p)$.\\

Let $E_\dde$ be the $\bb{Z}/p\bb{Z}$ torsor over $\padd$ deduced from $E$ by pulling back via the diagonal map. 

\begin {prop}
\label {formulaquadratic}
$\emph{a})$ If $E$ is skew symmetric, then there exists unique $r \in R$ such that $r(0) = 0$ and 
$$
r(x)^p - r(x) = xf(x)
$$

$\emph{b})$ The quadratic form $q: A \rightarrow \dis$ is then given by the formula $q(a) = r(a)$ for all $a \in A$.
\end {prop}

\begin {proof}
Indeed, since $E$ is skew-symmetric, $E_\dde$ is trivial, whence the existence of $r$ follows (uniqueness is obvious). Part (b) is clear (one argues as in the proof of prop.~\ref{formulabiadditive}/(b)).\qed
\end {proof}

\emph{Example}. We continue the previous example. Note that $xf(x) = x^{1+p} - x^{1+p^{-1}}$, whence $r(x) = x^{1+p^{-1}}$. Thus $q: \bb{F}_{p^2} \rightarrow \bb{Z}/p\bb{Z}$ is given by $a \mapsto a^{p+1} = \text{Nm}_{\bb{F}_{p^2}/\bb{F}_p}(a)$. This verifies thm. 1 in this special case.

\section* {Appendix B: Metric Groups}

We rapidly recall the notion of metric groups from \cite{drinfeldetel}. For more details, see Appendix A, \emph{loc. cit}.

We fix an algebraically closed field $K$ of chatacteristic zero. Let $A$ be a finite abelian group. A \emph{quadratic form} on $A$ (cf. \emph{loc. cit.}, sec. 2.11.1) is a map $q: A \rightarrow K^\times$ such that

i) The  function $B(a,b) = \frac{q(a+b)}{q(a)q(b)}$ from $A \times A$ to $K^\times$ is biadditive.

ii) $q(-a) = q(a)$ for all $a \in A$.\\

\emph{Remarks}. 1) It is easy to see that (i) and (ii) implies $q(na) = q(a)^{n^2}$ for all $n \in \bb{Z}$.

2) Let $p$ be a prime number and $\dis \inject K^\times$ be an injective group homomorphism. If $A$ is a $p$-group (i.e., \card($A$) is a power of $p$) equipped with a quadratic form $q$, then $q$ admits a (unique) factorization:
$$
A \xrightarrow{q^\p} \dis \inject K^\times
$$
In this case, one can equivalently define a quadratic from on $A$ to be a function $q^\p: A \rightarrow \dis$ such that

i$^\p$) The function $B^\p(a,b) = q^\p(a+b) - q^\p(a) - q^\p(b)$ from $A \times A$ to $\dis$ is biadditive.

ii$^\p$) $q^\p(-a) = q^\p(a)$.

\noindent This equivalent definition is the one used in this article.\\

The group $A$ equipped with the form $q$ is called a \emph{pre-metric} group. A homomorphism:
$$
\vp: (A,q) \rightarrow (A^\p,q^\p)
$$
of pre-metric groups is a homomorphism $A \xrightarrow{\vp} A^\p$ of abelian groups such that $q(a) = q^\p(\vp(a))$.

If $q(a) = 1$ implies $a = 0$, one says that $A$ is \emph{anisotropic}. A subgroup $H$ of $A$ is called \emph{isotropic} if $q|_H = 1$. In this case, one has $H \subset H^{\perp}$, where $H^{\perp}$ is the orthogonal complement of $H$ with respect to the biadditive map $B: A \times A \rightarrow K^\times$.  An isotropic subgroup $H$ is said to be \emph{Lagrangian} if $H = H^{\perp}$.

A pre-metric group $A$ is  called a \emph{metric group} if $B$ is nondegenerate. Assume henceforth that $A$ is a metric group. Let $H$ be an isotropic subgroup. One says that the group $H^\perp/H$ equipped with the induced quadratic form is a \emph{subquotient} of $A$. If $H$ is a \emph{maximal} isotropic subgroup, one can check that the resulting (anisotropic) metric group is independent, up to isomorphism, of the choice of $H$ (cf. \emph{loc. cit.}, Lemma A.31). One says that two metric groups are equivalent if the corresponding anisotropic metric groups are isomorphic. Equivalently, two metric groups are equivalent if and only if they have isomorphic metric subquotients. This defines an equivalence relation on the class of all metric groups. Let $W$ denote the set of equivalence classes. Let $[A]$ be the class of $A$ in $W$. One equips $W$ with an additive structure by setting:
$$
[A_1] + [A_2] = [A_1\oplus A_2]
$$
It is easy to check that this operation is well defined and endows $W$ with the structure of an abelian group (the inverse of $(A,q)$ being given by $(A,1/q)$).

\begin {define}
$W$ is called the \emph{Witt group} of metric groups.
\end {define}

\emph{Remark}. For each prime $p$, let $W_p$ be the subgroup of $W$ consisting of classes of metric $p$-groups. Then one has a decomposition:
$$
W = \bigoplus_p W_p
$$
It can be shown that $\card(W_2) = 16$ and $\card(W_p) = 4$ if $p$ is odd. For the proofs, see \emph{loc. cit.}, Prop. A.6 (for $p \neq 2$), Prop. A.13 (for $p=2$) as well as the remark following Def. A.33.


\begin{thebibliography}{ASM}
  
  
  \bibitem[B\'eg]{begueri} B\'egueri, L. {\em Dualit\'e sur un corps local a corps r\'esiduel alg\'ebriquement clos}, M\'em. Soc. Math. France (N.S.), 1980/81, no. 4.
  
  \bibitem[Boy] {mitya} Boyarchenko, M. {\em Characters of unipotent groups over finite fields}, Preprint, arXiv: \texttt{math.RT/0712.2614v1}, 2007.


  \bibitem [BD1] {motivated} Boyarchenko, M., Drinfeld V. {\em A motivated introduction to character sheaves and the orbit method for unipotent groups in positive characteristic}, Preprint, arXiv: \texttt{math.RT/0609769}, 2006.

  \bibitem [BD2] {foundations} Boyarchenko, M., Drinfeld V. {\em Character sheaves on unipotent groups in positive characteristic: foundations}, Preprint, arXiv: \texttt{math.RT/0810.0794v1}, 2008.

  \bibitem [Br1] {breen} Breen, L. {\em Biextensions altern\'ees}, Compositio Mathematica, {\bf 63}, pp. 99-122, 1987.

  \bibitem [Br2] {breen2} Breen, L. {\em Monoidal categories and multiextensions}, Compositio Mathematica, {\bf 117}, pp. 295-335, 1999.

  \bibitem[Des] {tanmay} Deshpande, T. {\em Heisenberg idempotents on unipotent groups}, in preparation.
  
  \bibitem[Dr]{drinfeld} Drinfeld, V. {\em Lectures on the geometric aspects of character theory for unipotent groups} (unpublished), University of Chicago, July 2005.
  
  \bibitem[DGNO] {drinfeldetel} Drinfeld, V., Gelaki, S., Nikshych, D., Ostrik, V. \emph{On braided fusion cateogories I}, to appear.
  
  \bibitem[Mil] {milne} Milne, J. {\em \'Etale cohomology}, Princeton University Series, \textbf{33}, Princeton University Press, 1980.
  
  \bibitem[Mum] {mumford} Mumford, D. {\em Bi-extensions of formal groups}, in \lq\lq Algebraic Geometry (Proceedings of Internat. Coll. Bombay, 1968)\rq\rq, pp. 307-322, Oxford University Press, 1969
  
  \bibitem[Ser] {serre} Serre, J-P. {\em Groupes Proalg\'ebriques}, Publ. Math. IHES, {\bfseries 7}, 1960.
  
  \bibitem[SGA7 I]{sga} Grothendieck, A. {\em S\'eminaire de g\'eom\'etrie alg\'ebrique du Bois-marie, 1967-1969, SGA7 I}, Lecture Notes in Math., {\bfseries 288}, Springer-Verlag, 1972.

  \bibitem[Tur] {turaev} Turaev, V. {\em Quantum invariants of knots and \emph{3}-manifolds}, W. de Gruyter, 1994.
\end{thebibliography}
\end {document}